\documentstyle[12pt]{amsart}

\begin{document}

\baselineskip=17pt

\title[]{Orbit equivalence and topological conjugacy of affine actions on
compact abelian groups}
\author[]{Siddhartha Bhattacharya}

\address{School of Mathematics, T.I.F.R, Mumbai 400005, India}

\email{siddhart@@math.tifr.res.in} 



\date{}
\maketitle 
\section{Introduction}
For a topological group $\Gamma $, by a $\Gamma $-flow we mean pair
$( X , \rho)$, where
$X$ is a topological space and $\rho $ is a continuous action of
$\Gamma $ on $X$. For any two $\Gamma $-flows $( X , \rho )$
and $( X^{'} ,
\sigma )$, a continuous map $f:X\rightarrow X^{'}$ is said to be
{\it $\Gamma$-equivariant \/} if 
$ f\circ \rho (\gamma ) = \sigma (\gamma )\circ f, \ \ \forall
\gamma
\in{\Gamma }$.
Two $\Gamma $-flows $( X , \rho )$     
and $( X^{'} ,
\sigma )$ are
said
to be {\it topologically conjugate \/} if there
exists a $\Gamma$-equivariant homeomorphism 
$ f : X \rightarrow X^{'} $ and
they are said to be {\it orbit equivalent \/} if there exists a
homeomorphism
$ f : X \rightarrow X^{'} $ which takes orbits under $\rho $ to orbits
under $\sigma $.

When $X , X^{'}$ are topological groups, an {\it affine \/} map from $X$
to
$X^{'}$ is a map of the form $ x \mapsto c\ \theta (x) $, where
$c\in{X^{'}}$ and $\theta $ is a continuous homomorphism from $X$ to
$X^{'}$. A $\Gamma $-flow $( X , \rho )$ is said to be affine if 
for all $\gamma $ in $\Gamma $, $\rho (\gamma )$ is an affine map. $(X ,
\rho )$
is said
to be an {\it automorphism flow \/}(resp. a {\it translation flow \/})
 if each 
$\rho (\gamma )$, $\gamma \in{\Gamma }$, is an automorphism on $X$ (resp
a
translation on $X$).
 $( X ,\rho )$ and $( X^{'} ,\sigma ) $
are said to be {\it algebraically conjugate \/} if there exists a
continuous
isomorphism $\theta : X \rightarrow X^{'} $ such that
$\theta \circ \rho (\gamma ) = \sigma(\gamma )\circ \theta $
$\ \forall \gamma \in{\Gamma }$.
\newline
\newline 
In this note we prove certain results concerning classification
of affine flows on compact connected metrizable abelian groups,
upto orbit equivalence and topological conjugacy.

We will denote by $S^{1}$ the usual circle group. For any locally compact
abelian group $G$, we denote
by $\widehat G$ the dual group of $G$. 
For a compact connected metrizable abelian group $G$,
we denote by
$L(G)$ the topological vector space consisting of all homomorphisms
from $\widehat G$ to $\mathbb{R}$, under pointwise addition and
scalar multiplication and the
topology of
pointwise convergence.
We define a map $E$ from $L(G)$ to $G$ by the condition
$$(\phi \circ E )(p) = e^{2\pi i p(\phi )} \ \ \ \forall p\in{L(G)} ,
 \phi \in{\widehat G}.$$
From the defining equation it is easy to see that $E$ is a continuous
homomorphism from
$L(G)$ to $G$. The kernel of $E$ can be identified with the set of all
homomophisms from $\widehat G$ to $\mathbb {Z}$; it is a totally
disconnected
subgroup of $L(G)$.
Note that using the duality theorem we can realise $L(G)$ with
set of all one-parameter subgroups of $G$ and $E$ with the map
$\alpha \mapsto \alpha (1)$. In particular when $G$ is a
 torus, $L(G)$   
can be identified with ${\mathbb R}^{n}$, the Lie
algebra of $G$, and $E$ can be identified with the standard exponential
map. However, in general $E$ is not
surjective, e.g for $G = {\widehat {\mathbb Q}}$,
where ${\mathbb Q}$ is the group of rational numbers equipped with
discrete topology, $L(G)$ is isomorphic to ${\mathbb R}$ and
consequently $E$ is not surjective.
\newline
\newline
Now let $\Gamma $ be a discrete group and $(G , \rho )$ be an affine
$\Gamma $-flow on $G$. Note that $\rho $ induces an automorphism
flow $\rho_{a} $ and a map $\rho_{t}:\Gamma \rightarrow G $
defined by
$$ \rho(\gamma )(x) = \rho_{a}(\gamma )(x)\ \rho_{t}(\gamma )
\ \ \forall x \in{G}.$$
  We define an automorphism action
$\rho_{*} $ of $\Gamma $ on $L(G)$ by
$$\rho_{*}(\gamma )(p)(\phi ) = p(\phi\circ\rho_{a}(\gamma ))
\ \forall\phi\in{\widehat G},\gamma\in{\Gamma}.$$
 In [ 1 ] and
[ 2 ] it was proved that any topological conjugacy between two
ergodic automorphisms of $n$-torus is affine. In [ 8 ] 
this was
generalized to certain
 class of affine transformations on compact connected
metrizable abelian groups. Here we prove the
following
\newline
\newline
{\bf Theorem 1 : }   
{\it Let $\Gamma $ be a discrete group and $G$ and $H$ be compact
connected
metrizable abelian groups. Let $\rho $, $ \sigma $ be affine actions
of $\Gamma $ on $G$ and $H$ respectively. Let $ f : G\rightarrow H$
 be a $\Gamma$-equivariant continuous map.
 Then there exist $c\in{H},$ a continuous homomorphism
$\theta : G \rightarrow H$ and a continuous map
$ S : G \rightarrow L(H)$ such that 
\newline

a) $S(e) = 0 $ and for all $x$ in $G$, the orbit of $S(x)$ under
$\sigma_{*}$
is bounded.

b) $ f(x) = c\  \theta (x) ( E\circ S ) (x),\ \ \forall x\in{G}$.
\newline
\newline
Moreover if $\rho $ and $\sigma $ are automorphism actions then 
$ S$ is a $\Gamma$-equivariant map from $(G,\rho)$ to
$(L(H),\sigma_{*})$.
 \/}
\newline
\newline
We show that under various additional conditions $S$ can be concluded
to be identically $0$, which means that $f$ is affine. This will be shown
to be the case, for instance, if the $\Gamma$-action on $H$ is
expansive (see corollary 1), see also Corollary 2 and
Remark 1 for other
applications of the theorem.
\newline
\newline 
In section 3 we classify translation flows upto orbit equivalence 
and topological conjugacy (in the continuous and discrete 
parameter cases respectively).
Firstly, for one-parameter flows of translations, generalizing a classical
result in the case of tori (see [ 3 ]) we prove the
following
\newline
\newline
{\bf Theorem 2 : }
{\it Let $G$ and $H$ be compact connected metrizable abelian groups and
$\alpha $ and $\beta $ be one-parameter subgroups of $G$ and $H$
respectively. Then the translation flows on $G$ and $H$ induced
by $\alpha $ and $\beta $ respectively are
orbit equivalent if and only if there exist a continous isomorphism
$\theta : G \rightarrow H$ and a nonzero $c\in{R} $ such that
 $\theta \alpha (t) = \beta (ct) \ \ \ \forall t\in{R} $. \/}
\newline 
\newline
We will also prove the following result which previously seems to have
been noted
only for ergodic translations (see [ 2 ]).
\newline
\newline
{\bf Theorem 3 : }
{\it Let $G$, $H$ be two compact connected metrizable abelian
groups and $\rho ,\sigma $ be two translation flows of a discrete   
group $\Gamma $ on $G$ and $H$ respectively. Then
$( G,\rho )$ and $( H,\sigma  ) $ are topologically
conjugate if and only if they are algebraically conjugate.   
\/}
\section{Rigidity of affine actions}
In this section 
we freely use various results from duality theory of locally compact
abelian groups; the reader is referred to [ 5 ] for details.
We will also use the following result due to VanKampen; for a proof see
[ 2 ], [ 7 ].
\newline
\newline 
{\bf Theorem }(VanKampen) {\it Let $G$ be a compact connected metrizable
abelian group and $f :G \rightarrow S^{1} $ be a continuous map.
Then there exist $ c \in{S^{1}}, $  $\phi \in{\widehat G}$  and a
continuous map
$h : G\rightarrow R $ such that
$$  h(0) = 0\ \ ,\ \ f(x) = c\ \phi (x)\ e^{2\pi ih (x)} \ \ \forall
x\in{G}.$$
Moreover $c,\phi $ and $h$ are uniquely defined. \/}
\newline
\newline
Now for any two groups $G , H$ and any continuous map 
$ f : G \rightarrow  H $ we will define a continuous
homomorphism $ \theta (f) : G \rightarrow  H $ as follows.
For each character $\phi $ of $H$ let
$c_{\phi }\in{S^{1}},$
${\widehat \theta}(\phi )\in{\widehat G}$
and $f_{\phi} : G\rightarrow R $ be such that
$$ f_{\phi} (0) = 0 \ ,\
 \phi \circ f (x) = c_{\phi }\ {\widehat \theta}(\phi )(x)\ e^{2\pi i
f_{\phi
}(x)}\ \
\forall
x\in{G}; $$ note that by VanKampen's theorem there exist 
$c_{\phi },{\widehat \theta}(\phi)$
and $f_{\phi }$ satisfying the conditions and they are
  unique. From the uniqueness one can deduce that
$ \phi \mapsto {\widehat \theta}(\phi )$ is a homomorphism
 from $\widehat H$ to $\widehat G$. By the duality theorem there
exists a  continuous homomorphism  
$\theta (f): G \rightarrow H$ such that
$${\widehat \theta}(\phi ) = \phi \circ \theta (f)\ \ \ \forall \phi \in 
{\widehat H}.$$ 
Using the uniqueness part of VanKampen's theorem it is easy to 
see that 
\newline

i) If $f$ is a continuous homomorphism then $ \theta (f) = f $.

ii) If $ f : G_{1} \rightarrow G_{2}$ and
$ g : G_{2} \rightarrow G_{3}$
be two continuous maps then
 
$ \theta ( g\circ f ) = \theta( g)\circ \theta (f) $.
\newline

The following corollary of VanKampen's theorem is essentially due to Arov
(see [ 2 ]) .
\newline
\newline
{\bf Proposition 1 : }
{\it Let $\Gamma $ be a discrete group and $ \rho : \Gamma\rightarrow
Aut(G)$
and $ \sigma : \Gamma\rightarrow Aut(H)$
be  automorphism actions  of $\Gamma $ on $G$ and $H$ respectively.
 Then $\rho $ and $\sigma $ are topologically conjugate if and only
if they are algebraically conjugate. \/}
\newline
\newline
{\bf Proof :} Let $f$ be a topological conjugacy between $\rho $ and
$\sigma $.
Using i) and ii) we see that 
$$ \theta (f) \circ \rho (\gamma ) = \sigma (\gamma )\circ \theta (f) 
 \ \forall \gamma\in{\Gamma }.$$
Since $f$ is a homeomorphism, it follows that 
$\theta (f) $ is an isomorphism. Hence $\rho $
and $\sigma $ are algebraically conjugate.
\newline
\newline
The following lemma generalizes of VanKampen's theorem.
\newline
\newline
{\bf Lemma 1 : }
{\it Let $G$, $H$ be two compact connected metrizable abelian groups
and $f$ be a  continuous map from $G$ to $H$.
 Then there exist $c\in{H},$ a continuous homomorphism 
$\theta : G \rightarrow H$ and a continuous map 
$ S : G \rightarrow L(H)$ such that $S(0) = 0 $ and
$ f(x) = c\  \theta (x) ( E\circ S ) (x),\ \ \forall x\in{G}$.
Moreover $c, \theta $ and $S$ are unique. \/}
\newline
\newline
{\bf Proof : }
For each character $\phi $ of $H$ define 
$c_{\phi }\in{S^{1}},\theta ^{'}(\phi )\in{\widehat G}$
and $f_{\phi} : G\rightarrow R $ by the condition
$$ f_{\phi} (0) = 0 \ ,\
 \phi \circ f (x) = c_{\phi }\ \theta^{'}(\phi )(x)\ e^{2\pi i f_{\phi
}(x)}\ \
\forall 
x\in{G};$$ note that by VanKampen's theorem there exist uniquely
defined 
$c_{\phi }, \theta^{'}(\phi )$ and 
$f_{\phi }$ satisfying the condition.
From the uniqueness it follows that $$f_{\phi \psi} =
f_{\phi} + f_{ \psi}\ ,\
f_{\phi^{-1}} = -f_{\phi}.$$ Define $ S : G \rightarrow L(H)$ 
by $ S(x)(\phi ) = f_{\phi }(x) $. Since each $f_{\phi }$ is 
continuous, $S$ is continuous. Similarly using uniqueness
of $c_{\phi }$  we see that the map 
$\phi \mapsto  c_{\phi }$ is a homomorphism
from $\widehat H$ to $ S^{1}$.
By the duality theorem there exists $c\in{H}$ 
such that $c_{\phi } = \phi (c) \ \ \forall \phi\in{\widehat H}$.
Also putting $\theta = \theta (f)$ we see that
 $\theta^{'}(\phi ) = \phi \circ \theta , $
$\ \ \forall \phi \in{\widehat H}$. Hence
for all
$x\in{G} $ and
$ \phi \in{\widehat H}$, we have
$$ \phi \circ f (x) = c_{\phi }\ \theta^{'}(\phi ) (x)\ e^{2\pi i f_{\phi
} (x)}
 = \phi (c)\ (\phi \circ \theta ) (x)\ ( \phi \circ E \circ S ) (x).$$
Since characters separate points, 
$ f(x) =c\ \theta (x) ( E\circ S ) (x),\ \ \ \forall x\in{G}$.
Using VanKampen's theorem we see that for a fixed 
$\phi \in{\widehat H}$, $\phi (c),\phi \circ \theta  $
and the map $ x\mapsto S(x)(\phi ) $ are determined by the equation
$$ (\phi \circ f) (x) = \phi (c)\ (\phi \circ \theta)  (x)\
e^{2\pi i
S(x)(\phi )}.$$
Hence $c,\theta $ and $S$ are unique.
\newline
\newline
{\bf Proof of Theorem 1 : }
    Suppose $ f = c\ \theta (E\circ S )$ where $ c, \theta $ and
$ S$ are as in Lemma 1.  Fix any $\gamma \in{\Gamma}$. 
Note that for all
$x \in{G}$,
$$ f\circ \rho(\gamma ) (x)
= c_{1 }\ \theta_{1}(x) (E\circ S_{1})(x),$$
where
$\ c_{1} = f\circ \rho (\gamma )(e)$,
 $\theta_{1} = \theta \circ \rho_{a}(\gamma )$ and
$S_{1}(x) = S\circ \rho(\gamma ) (x) - S\circ\rho(\gamma ) (e)$.
Also for all $x\in{G}$,
$$ \sigma (\gamma )\circ f(x)
= c_{2 }\ \theta_{2}(x) (E\circ S_{2 })(x),$$
where
$c_{2} = \sigma (\gamma )\circ f (e)$
, $\theta_{2} = \sigma_{a}(\gamma )\circ \theta $
 and $S_{2 } = \sigma_{*}(\gamma )\circ S $.
From the uniqueness part of Lemma 1 it follows that
$S_{1} = S_{2} $ i.e.
$$ S \circ \rho (\gamma )(x) - S \circ \rho(\gamma )(e) =
 \sigma_{*}(\gamma )\circ S (x),\ \ \forall x\in{G}.$$
Since for a fixed $x\in{G}$ the left hand side is contained
 in a bounded subset of $L(H)$,
it follows that for all $x$ in $G$, the $\sigma_{*}$-orbit 
 of $S(x)$ is bounded. Also it is easy to see from the previous
identity that when $\rho $ and $\sigma $ are automorphism actions, $S$
is a $\Gamma$-equivariant map.
\newline
\newline
When $X$ is a topological group,
$( X ,\rho ) $
 is said to be {\it expansive \/} if there exists a neighbourhood $U$ of
the identity
 such that for any two distinct elements $x,y\in{X}$ there exists
a $\gamma \in{\Gamma }$
such that $\rho (\gamma )(x)\ \rho (\gamma ) (y)^{-1} $ is not
contained in $U$; such a neighbourhood is called an expansive
neighbourhood. For various characterizations of expansiveness
of automorphism actions on compact abelian groups the reader
is referred to [ 6 ].  
\newline
\newline
{\bf Corollary 1 :}
{\it Let $\Gamma $ be a discrete group and $G$ and $H$ be compact
connected
metrizable abelian groups. Let $\rho $, $ \sigma $ be affine actions
of $\Gamma $ on $G$ and $H$ respectively such that 
$( H, \sigma )$ is expansive. Then every $\Gamma$-equivariant
continuous map  $ f :(G, \rho )\rightarrow (H , \sigma )$
 is an affine map. \/}
\newline
\newline
{\bf Proof : }
Since $\sigma $ is expansive, $\sigma_{a}$ is an expansive
automorphism action on $H$. We claim that for every nonzero point
$p\in{L(H)}$, the orbit of $p$ under $\sigma_{*}$ is unbounded. Suppose
not.
Choose a non-zero $p$ and a compact set $C\subset L(H)$
such that orbit of $p$ under $\sigma_{*}$ is contained in $C$.  
Since kernel of
$E$
is totally disconnected, there exists a sequence $\{ t_{i}\}$ such
that $t_{i}\rightarrow 0 $ as $i\rightarrow \infty $ and
$E(t_{i}p)\ne e \ \forall i $. Let $U$ be an expansive neighbourhood
of $e$ in $H$. Since $e$ is fixed by $\sigma_{a} $, $e$ is the only
element
in $G$ whose orbit under $\sigma_{a} $ is contained in $U$.
Since $t_{i}\rightarrow 0 $, it is easy
to see that $ \cup\ t_{i}^{-1} E^{-1}(U) = L(H)$. From the compactness
of $C$ it follows that there exists $n$ such that
$t_{n}C \subset E^{-1}(U)$ i.e. $E(t_{n}C) \subset U $.
Since the orbit of $t_{n}p$ under $\sigma_{*}$ is contained in $t_{n}C$
and
$ E\circ \sigma_{*}(\gamma ) = \sigma_{a} (\gamma )\circ E ,\ \forall
\gamma$,
this implies that orbit of $E(t_{n}p )$ under $\sigma_{*} $ is contained
in $U$. 
This contradicts the fact that $E(t_{n}p)\ne e$.

Now suppose $f = c\ \theta (E\circ S)$, where $c,\theta $ and $S$ are
as in Theorem 1. From Theorem 1 and the previous argument it follows
that $S = 0$ i.e. $f$ is an affine map.
\newline
\newline
For a set $A$ we denote by $|A|$ the cardinality of $A$.
\newline
\newline
{\bf Corollary 2 :}
{\it Let $\Gamma $ be a discrete group and $G$ and $H$ be compact
connected
metrizable abelian groups. Let $\rho $, $ \sigma $ be 
automorphism actions
of $\Gamma $ on $G$ and $H$ respectively such that

a ) $\{g\in{G}\ |\ |\rho(\Gamma )(g)|<\infty\}$ is dense in $G$.

b ) for any natural number $k$, the set
 $\{h\in{H}\ |\ |\sigma(\Gamma)(h)| = k \}$ 
is

\ \ \ \  totally
disconnected.
 \newline
Then every continuous $\Gamma$-equivariant map 
$f : ( G,\rho )\rightarrow ( H,\sigma )$ is an affine map. \/}
\newline
\newline
{\bf Proof :}
Suppose $f = c\ \theta (E\circ S)$, where $c,\theta $ and $S$ are
as in Theorem 1. Let $g\in{G}$ be such that the orbit of $g$
under $\rho $ is finite. Since $S$ is
$\Gamma$-equivariant by Theorem 1, the orbit of $S(g)$ under 
$\sigma_{*}$ is also finite. Since $\sigma_{*}$ is a linear 
action on $L(H)$ and $E$ is a $\Gamma$-equivariant map from
$(L(H),\sigma_{*})$ to $(H,\sigma)$, the orbit of 
$E(t S (g))$ under $\sigma $ is finite for all $t\in{\mathbb R}$.
Now from  b) it follows that $S(g) = 0$. Since 
$\{g\in{G}\ |\ |\rho(\Gamma)(g)|<\infty\}$ is dense in $G$,
 $S = 0$ i.e. $f$ is an affine map.
\newline
\newline
{\bf Remark 1 : }
It is easy to see that condition (a) as in corollary 2,
holds when $G = T^{n}$ for some $n$.
Various other conditions  under which the
set of periodic
orbits of an automorphism action on a compact abelian group is dense, viz
condition (a) as in corollary 2 holds, are described in [ 4 ].
Condition (b) holds in the case of $H = T^{n}$ if $\Gamma $
contains an element acting ergodically; more generally this
holds for any finite dimensional compact abelian group $H$. 
\section{Classification of translation flows} 
{\bf Lemma 2 : }
{\it Let $G$ be a compact connected metrizable abelian group and
$\alpha $ be a one-parameter subgroup of $G$. Then there exists
a $p\in{L(G)}$ such that 
$E(tp) = \alpha (t)\ \forall t\in{\mathbb R}$.\/}
\newline
\newline
{\bf Proof : }
For each $\phi \in{\widehat G}$, we define
 $\alpha_{\phi }\in{\mathbb R}$ by
$$ \phi\circ\alpha (t) = e^{2\pi i \alpha_{\phi }t} \ \forall
t\in{\mathbb R}.$$
Since $\phi\circ \alpha $ is a continuous homomorphism from 
$\mathbb R$ to $S^{1}$, $\alpha_{\phi }$ is well defined.
We define 
$p\in{L(G)}$ by 
$p(\phi ) = \alpha_{\phi }\ \forall\phi\in{\widehat G}$. Fix
any $t\in{\mathbb R}$. From the defining equation of $E$ it follows that
for all $ \phi\in{\widehat G}$,
$$\phi\circ E(tp) = e^{2\pi i t p(\phi )} = e^{2\pi i \alpha_{\phi }t}
= \phi\circ\alpha (t).$$
Since characters separate points, it follows that 
$\alpha (t) = E(tp)\ \forall t\in{\mathbb R}$.
\newline
\newline
 The following
lemma is needed to prove Theorem 2. The main idea of
the Proof is derived from [ 3 ].
\newline
\newline
{\bf Lemma 3 : }
{\it Let G be a compact connected metrizable abelian group and
$p,q\in{L(G)},$ with $p\ne 0$.
Let $f : R\rightarrow L(G) $ be a bounded continuous function such that
 $$f(0) = 0\ ,\ \{E(\ tp+f(t)\ ) \ |\ \ t\in{R}\} = \{E(\ tq\ )\ |\
\ t\in{R}\}.$$
Then $p = cq$ for some nonzero $c\in{R}$. \/}
\newline
\newline
{\bf Proof : }
First we will prove the special case when $G = T^{2}$, the
two-dimensional torus. After suitable
identifications we have 
$$L(G) = R^{2}\ ,\  
E(x_{1},x_{2}) = Exp(x_{1},x_{2}) = 
( e^{2\pi i x_{1}},e^{2\pi i x_{2}}).$$
Define a  function $d : R^{2} \rightarrow R^{+}$ by

$d(x) = $ distance between the point x and the 
line $\{tq\ |\ t\in{R}\}$

$ \ \ \ \ \ \ $ = inf $\{||y||\ |\ x+y =tq $ for some $ t\in{R}\}$.
\newline
By our hypothesis , for all $ t\in{R}$ , $tp + f(t) = z + t^{'}q $
for some $t^{'}\in{R} , z\in{{\mathbb{Z}}^{2}}$.
This implies
$$\{d(\ tp+f(t)\ )\ |\ t\in{R}\}\subset \{d( z )\ |\
z\in{{\mathbb{Z}}^{2}}\}$$ 
Since the map $ t\mapsto d(tp+f(t))$ is continuous, the left hand
side is 
a connected subset of $\mathbb{R}$ containing $0$. Since the right hand
side is
countable, $ d(tp+f(t)) = 0, \ \ \forall  t $.
Since $f$ is bounded this implies that $d(tp)$ is bounded by a 
constant $M$, for all $t\in{R}$.
Since distinct lines in $R^{2}$ diverge from each other
 we conclude that $p = cq$ for some $c\ne 0 $.
\newline

To prove the general case choose $\phi $ such that $p(\phi ) \ne 0 $.
For
each $\psi \in{\widehat G}$ define
$h : G \rightarrow T^{2} $ and $h^{*} : L(G) \rightarrow R^{2} $
by 
$$ h(x) = (\phi (x)\ ,\ \psi (x))\ ,\ \  
h^{*}(r) = (r(\phi)\ ,\ r(\psi )).$$
Now for all $r\in{L(G)}$,
$$h\circ E(r) = (\phi\circ E(r),\psi\circ E(r))\ ,\ 
Exp\circ h^{*}(r) = (e^{2\pi i r(\phi )},e^{2\pi i r(\psi )}).$$
From the defining equation of $E$ it follows that 
$h\circ E = Exp\circ h^{*}$. Define $p^{'},q^{'}\in{ R^{2}}$ and
$f^{'} : G\rightarrow R^{2}$ by 
$p^{'} = h^{*}(p)$, $ q^{'} = h^{*}(q)$ and
$f^{'} = h^{*}\circ f$.
From our hypothesis it follows that
$$ \{Exp(tp^{'}+f^{'}(t))\ |\ t\in{R}\}
= \{Exp(tq^{'})\ |\ t\in{R}\}$$
 Now applying the special case we see that
$(p(\phi ),p(\psi )) =  c\ (q(\phi ),q(\psi ))  $ for some nonzero
real number $c$. Therefore $q(\phi ) \ne 0 $ and, since $\psi $ is
arbitrary,
 $ p = c_{0}q\ $
 where $c_{0} = p(\phi ) /q(\phi ) $. 
\newline
\newline
{\bf Proof of Theorem 2 : }
Let $h$ be an orbit equivalence between the translation flows induced by 
$\alpha $ and $\beta $. Define $f:G\rightarrow H$ by
$f(x) = h(x)h(e)^{-1}$. Then $f(e) = e $ and it is easy to check that $f$
is also an orbit equivalence. Suppose $ f = \theta ( E \circ S )$, where 
$\theta $ and $S$ are as in Lemma 1. By Lemma 2  there exists 
$p ,q$ in $L(H)$ such that 
$ E(tq) = \beta (t) $ and $ E(tp) = \theta (\alpha (t))$.
Since $f$ is
an orbit equivalence,
$$\{f(\alpha (t))\ |\ t\in{R}\} = \{\beta (t)\ |\ t\in{R}\}$$
Now for all $ t \in{R}$, $\beta (t) = E(tq) $ and
$$ f(\alpha (t)) = 
(\theta \circ \alpha) (t) (E\circ S \circ \alpha )(t) 
= E(tp +S\circ \alpha (t)).$$

Since $f$ is a homeomorphism, $\theta $ is an isomorphism.
Hence $\theta (\alpha ) \ne 0 $, i.e. $p \ne 0 $. By applying Lemma 3
 we see
that  $ p = cq $ for some nonzero $c$. Therefore for some
$ c\ne 0 $,
 $\theta ( \alpha (t)) = \beta (ct), \ \ \ \forall t \in{R}$.   
This proves the theorem.
\newline
\newline
{\bf Proof of Theorem 3 : }
Let $f$ be a topological conjugacy between the induced translation
flows. Suppose $ f = c\ \theta \ ( E \circ S ) $, where
$c,\theta$ and $ S $ are as in Lemma 1. Since $f$ is a
homeomorphism, $\theta $ is an isomorphism. We claim that
$\theta $ is an algebraic conjugacy between
$( G,\rho )$ and $( H,\sigma ) $. To see this fix any      
$ \gamma \in{\Gamma}$. Define $x_{0} \in{G}$ ,
$ y_{0}\in{H} $ by
$$ \rho (\gamma )(x) = x_{0} x\ \ ,\ \ \sigma(\gamma )(y)
= y_{0} y \ \ \forall x\in{G},\ y\in{H}.$$
Then for all $x\in{G}$,
$$ f(x_{0} x )
= c\ \theta ( x_{0} x ) (E\circ S)(x_{0} x )
= c_{1}\ \theta (x) (E\circ S_{1} )(x), $$
where
$ c_{1} = c\ \theta (x_{0} )(E\circ S)$
$(x_{0} )\ ,$
$\ S_{1}(x) = S(x_{0} x) - S(x_{0} )$.
Also for all $x\in{G}$,
$$ y_{0} f( x )
= c_{2}\ \theta (x) (E\circ S)(x),$$
where $c_{2} = c\ \theta (x_{0} )$.
From the uniqueness part of Lemma 1 it follows that $S_{1} = S $
i.e.
$$ S (x_{0} x ) = S(x) + S (x_{0} )
,\ \ \forall x\in{G}. $$
Putting $ x = x_{0} , x_{0}^{2} ,\dots$
and using the above recursion relation we obtain
$$ S (x_{0}^{n}) = n S(x_{0}), \ \ \forall n\in{{\mathbb{Z^{+}}}}.$$
Since the left hand side is contained in the image of $S$, which
is compact, it follows that $ S(x_{0} ) = 0 $. Since
$ c_{1} = c_{2}$, this implies
$ \theta ( x_{0} ) = y_{0}$. Hence
$\theta \circ \rho(\gamma ) = \sigma(\gamma )\circ \theta $.
\newline
\newline
\newline
{\bf Acknowledgement : } I wish to thank Prof. S.G Dani for his
helpful suggestions and valuable help during preparation of this note. 
\newline
\newline
{\bf References :}
\newline
\newline
[ 1 ]\ \ R.L Adler, R. Palais (1965) Homeomorphic conjugacy
of automorphisms on the torus. Proc. Amer. Math Soc.
{\bf 16}: 1222-1225
\newline
[ 2 ]\ \ D.Z Arov (1963) Topological similitude of automorphisms and
translations on compact abelian groups.
 Uspehi Mat.Nauk . {\bf 18}: 133-138
\newline
[ 3 ]\ \  D. Benardete (1988) Topological equivalence of flows on
homogeneous
spaces,
and divergence of one-parameter subgroups of Lie groups.
 Trans.Amer.Math. Soc.  {\bf 306}: 499-527 
\newline
[ 4 ]\ \ B. Kitchens, K. Schmidt (1989) Automorphisms of compact
groups.
Ergodic Theory Dynamical Systems. {\bf 9}: 691-735 
\newline
[ 5 ]\ \ S.A Morris (1977) Pontryagin duality and the structure  
   of locally compact abelian groups.
 London Mathematical Society Lecture
   Note Series {\bf 29}: Cambridge University Press
\newline
[ 6 ]\ \ K. Schmidt (1990) Automorphisms of compact abelian groups and
affine varieties.  Proc. London Math. Soc.  {\bf 61}: 480-496   
\newline
[ 7 ]\ \  E.R VanKampen (1937) On almost periodic functions of constant
absolute value.
 J.London.Math. Soc. {\bf 12}: 3-6 
\newline
[ 8 ]\ \  P. Walters (1969) Topological conjugacy of affine
transformations on
compact abelian groups.
Trans.Amer.Math. Soc. {\bf 140}: 95-107
\end{document}